# All-to-all Routing on Vertex Transitive Digraph Networks

Nyumbu Chishwashwa[1], Vance Faber[2,*] and Noah Streib[3]

[1] University of Western Cape, Cape Town, South Africa

[2] Center for Computing Sciences, Bowie, Maryland

[3] Center for Computing Sciences, Bowie, Maryland

[*]vance.faber@gmail.com

**Abstract.** We discuss an open problem for all-to-all routing on highly symmetric digraph networks; namely, which digraphs have a spanning factorization.

## Introduction

In this paper, we define a *regular digraph* $D$ to be a set of vertices $V$ and a set of directed edges $E$ such that the in-degree and out-degree of every vertex is a constant $d$ and there is a directed path between any two vertices. These digraphs are often employed as the model for network connections in parallel computing, and a critical problem in that regard is optimization of the network with respect to global communication tasks, including all-to-all routing. See for example [3]. In this note, we utilize factorizations of regular digraphs developed in [1] to extend routing schemes used in highly symmetric digraphs to a wider class. We discuss the question posed in [1] whether every vertex transitive digraph has a spanning factorization, which would lead to a conflict-free all-to-all routing scheme for these graphs.

## Factorizations

A *factor* of a regular digraph $G$ is a function (written as applying from the right instead of the usual left) on the vertices which maps each vertex to a directed neighbor. We also think of a factor as a set edges of $G$. A factor is a *1-factor* if it is one to one. A *factorization* $F = \{F_1, F_2, \ldots, F_d\}$ of a regular digraph of degree $d$ is a disjoint decomposition of the edges into factors. A *1-factorization* is a factorization where all factors are 1-factors. In a factorization, each vertex has one out-edge in every $F_i$. In a 1-factorization, each vertex has one out-edge and one in-edge in every $F_i$. Each 1-factor $F_i$ is a disjoint union of directed cycles which contains all the vertices of the graph. By Petersen's theorem (see for example [4]), every regular digraph has a 1-factorization. For completeness, we provide a proof.

*Petersen's Theorem. Every directed graph $G$ where the in-degree and out-degree of every vertex is $d$ has an edge disjoint decomposition into $d$ 1-factors.*

*Proof.* Form an auxiliary graph $B$ with two new vertices $u'$ and $u''$ for each vertex $u$. The edges of $B$ are the pairs $\{u', v''\}$ where $(u, v)$ is a directed edge in $G$. The undirected graph $B$ is bipartite and regular with degree $d$ and so by Hall's Marriage Theorem, it can be decomposed into $d$ 1-factors. Each of these 1-factors corresponds to a directed 1-factor in $G$.

*Definitions*. Let $F_1$, $F_2$, $\cdots$, $F_d$ be the factors in a factorization of $G$. We call a finite string of symbols from the set $\{F_i\}$ a *word*. So given a vertex $u$ and factor $F_i$, $uF_i$ is a neighbor. Then each word $\omega$ is the composition of the 1-factors appearing in the word so that $u\omega$ is the vertex which is the result of following the directed path from $u$ given by $\omega$. It is convenient to write the function which is the composition of two words $\alpha$ and $\beta$ as $\alpha\beta$. Given a factorization and a vertex $v$ in $G$, we can easily find a set of words *generated* by the $F_i$ of the form

$$W = \{\omega_0 = \varnothing, \omega_1 = F_1, \omega_2 = F_2, \ldots, \omega_d = F_d, \omega_{d+1}, \ldots, \omega_{n-1}\}$$

such that the $v\omega_i$ include all the vertices of $G$. (Note that in this definition, we require the words $\omega_i$ for $1 \le i \le d$ to be the factors $F_i$. This requirement is absent in [1].) If for each vertex $v$ we have a set $W(v)$ generated by the $F_i$, the collection of $W(v)$ form an *all-to-all routing* on $G$. If there is a single set of words $W$ such that $W(v) = W$ for all vertices $v$ forms an all-to-all routing, then we say that the factorization (1-factorization) and the set of words $W$ form a *spanning factorization* (*1-factorization*) of $G$. Note that given a vertex $v$, a word $\omega$ in $W(v)$ and a factor $F_i$ in $\omega$, then $F_i$ represents a unique edge in $G$.

**Spanning 1-Factorizations**

This section reviews the definitions and major result from [1].

*Definition*. A *schedule* associated to a set of words $W$ *generated* by the $F_i$ is an assignment of a time (a label) to each occurrence of each 1-factor in the words of $W$ such that no time is assigned more than once to a particular 1-factor and times assigned to the 1-factors in a single word are increasing. The *makespan of a schedule* is the largest time assigned to any of the 1-factors.

The power of a spanning 1-factorization lies in the fact that a schedule for a single $W$ can be used to describe an algorithm for conflict free all-to-all routing. The schedule also provides an algorithm for other communication primitives like global broadcast or summation. This is discussed at length in [1] where the following result is proved.

*Spanning 1-Factorization Theorem. Suppose we have a schedule for a set of words W of a spanning 1-factorization of the graph $G$. Then the collection of directed label-increasing paths $v\omega_i$ for all $v$ and non-empty $\omega_i$ have the property that no edge in the graph is assigned the same time twice. A schedule for a spanning 1-factorization yields a time labeled directed path between every two vertices so that no edge is labeled with the same time twice.*

*Proof.* Each edge in the graph is assigned to a single 1-factor. Assume there is an edge in the 1-factor $F_i$ that has been assigned the same time twice. Since every occurrence of $F_i$ in the words in $W$ has been assigned a unique time, this can only mean that there are two different vertices $u$ and $v$ and an initial subword $\omega$ of a word in $W$ such that the edges $(u\omega, u\omega F_i)$ and $(v\omega, v\omega F_i)$ are the same edge. Then $u\omega$ and $v\omega$ must be the same vertex. Let us assume that this is the shortest $\omega$ for which this happens. The word $\omega$ cannot be empty since $u$ and $v$ are different. But then the last 1-factor in $\omega$ must contain two distinct edges with the same endpoint, contradicting the definition of a 1-factor. If we start with a spanning 1-factorization, then all the non-empty paths from $v$ are unique, there are $n-1$ of them and none of them can return to $v$ so they must reach to every other vertex in the graph.

**Vertex Transitive Digraphs**

A digraph $G$ is *vertex transitive* if for any two vertices $u$ and $v$ there is an automorphism of $G$ which maps $u$ to $v$.

**We examine the relationship between the properties spanning 1-factorization and vertex transitive.**

We use the following Cayley coset representation of a vertex transitive digraph (see [1]).

*Definition* (Cayley coset digraph). Let $\Gamma$ be a finite group, $\mathrm{H}$ a subgroup and $S$ a subset. Suppose

1) $S \cap H = \varnothing$,
2) $HSH \subseteq SH$,
3) $S$ is a subset of distinct coset representatives of $H$ in $\Gamma$.

Then we can form the *Cayley coset digraph* $G = (\Gamma, S, H)$ with the cosets $\{gH : g \in \Gamma\}$ as vertices and the set of pairs $g(H, sH) = (gH, gsH)$ with $s \in S$ as edges. When $H$ is the identity subgroup, the graph is a *Cayley digraph*.

A minor modification of the proof of Sabidussi [5] shows that a digraph is vertex transitive if and only if it is a Cayley coset digraph. For completeness, we provide a proof. An important aspect of the proof shows that one can construct a Cayley coset digraph from a vertex transitive digraph by using the automorphism group as the group $\Gamma$ required in the definition and the subgroup of automorphisms that fix a vertex as the required subgroup $H$. The generators $S$ correspond to automorphisms that map a vertex to a neighbor.

*Cayley Coset Digraph Characterization.* A graph is vertex transitive if and only if it is isomorphic to a Cayley coset digraph.

*Proof.* First, we show that a Cayley coset digraph is vertex transitive. Let $G = (\Gamma, S, H)$. We will show that $G$ is vertex transitive.

Suppose $f, g \in \Gamma$ and $fH = gH$.

P1: For every $s \in S$ there exists $t \in S$ such that $fsH = gtH$.

To prove this, note that $fH = gH$ implies $g^{-1}fH = H$ so there exists $h \in H$ with $g^{-1}f = h$. Because $HSH \subseteq SH$, there exists $t \in S$ such that $g^{-1}fsH = tH$ which shows $fsH = gtH$.

P2: For each $f \in \Gamma$ the function $\theta_f$ defined by $\theta_f(aH) = faH$ is an automorphism of $G$.

We have to show that if $bH = asH$ for $s \in S$ and $faH = cH$ then there exists $t \in S$ such that $\theta_f(bH) = ctH$. But $faH = cH$ implies that $aH = f^{-1}cH$ so by P1, $bH = asH = f^{-1}ctH$ for some $t \in S$. Thus $\theta_f(bH) = fbH = ctH$.

Clearly, P2 implies $G$ is vertex transitive because we can map any coset $aH$ to a coset $bH$ with the automorphism $\theta_{ba^{-1}}$.

Second, we show that a vertex transitive digraph is isomorphic to a Cayley coset digraph. To prove this, we show that a vertex transitive digraph $D$ is isomorphic to a graph $G = (\Gamma, S, \mathrm{H})$ where $\Gamma$ is the automorphism group of $D$, $S$ is a set of distinct coset representatives for the automorphisms which map a fixed vertex $v_0$ to its neighbors and $H$ is the set of automorphisms which fix $v_0$. To show that $G$ satisfies our definition of a Cayley coset digraph, we just need to verify that it satisfies 1 through 3. Only 2) is not obvious (but it is inherent in the classical proof.) We consider $s \in \Gamma$, $s(v_0) \in N(v_0)$, the neighbors of $v_0$. Let $h \in H$. Because $(v_0, s(v_0))$ is an edge and $h$ is an automorphism fixing $v_0$, then $(v_0, hs(v_0))$ is also an edge. This means there is some $t \in S$ such that $hsH = tH$. Consequently, $HSH \subseteq SH$. Let $\theta$ between $(\Gamma, S, \mathrm{H})$ and $D$ be given by $\theta(g(v_0)) = gH$. We show that $\theta$ is well-defined and is a graph isomorphism. If $g(v_0) = f(v_0)$ then $f^{-1}g(v_0) = v_0$ so $f^{-1}g \in H$ and thus $fH = gH$. Suppose $(v, w)$ is an edge in $D$ and let $\alpha$ be an automorphism such that $\alpha(v_0) = v$. The edges $(\alpha(v_0), \alpha s(v_0))$ with $s \in S$ are distinct since the edges $(v_0, s(v_0))$ are, so there is some $s \in S$ such that $\alpha s(v_0) = w$. Thus $(\theta(v), \theta(w)) = (\theta\alpha(v_0), \theta\alpha s(v_0)) = (\alpha H, \alpha sH)$ and $\theta$ maps edges to edges.

*Connected Case*. In addition, $G = (\Gamma, S, H)$ is connected if and only if

*4)* $\Gamma = \langle S \rangle H$.

*Proof.* It is obvious that if $\Gamma = \langle S \rangle H$, then $G$ is connected. Conversely, suppose $G = (\Gamma, S, H)$ is connected but $\Gamma \neq \langle S \rangle H$ and consider the induced subgraph $G_H$ generated by vertices $gH$ with $g \in \langle S \rangle H$. Then since $G$ is connected, there must be an $fH$ not in $G_H$ and $g \in \langle S \rangle$ such that $fH = gsH$ for some $s \in S$. This means that $f = gsh$ for some $h$ and so $f \in \langle S \rangle H$ after all.

*Example of what can go wrong*. Given a connected vertex transitive digraph we want to generate a spanning 1-factorization. It would be nice if each edge generator in $S$ created a 1-factor. But that may not be the case unless the coset representatives are carefully chosen. In [1], we showed how to find choices that yield a 1-factorization for several classes of vertex transitive digraphs. In this paper, we exhibit properties that allow us to do this for a more general set of vertex transitive graphs. Here, we give an example of bad choices. We start with a directed version of the Petersen graph with degree 2. We use the Cayley coset definition. The group is a semidirect product group on the two generators $\alpha$ and $\theta$ with

$$\alpha^4 = 1 = \theta^5$$
$$\alpha^{-1}\theta\alpha = \theta^2.$$

The subgroup is $H=\{1,\alpha^2\}$. The edges in the graph are the cycle $C=(H,\theta H,\theta^2 H,\theta^3 H,\theta^4 H)$, the edges $(\theta^i H,\theta^i\alpha H)$ and $(\theta^i\alpha H,\theta^i H)$ and the cycle $C'=(\alpha H,\alpha\theta H,\alpha\theta^2 H,\alpha\theta^3 H,\alpha\theta^4 H)$. (Note that it takes some calculation to determine the labels for this last set of cosets.) If you form a 1-factor with these coset representatives, the generators $\alpha$ and $\theta$ each form a 1-factor of the graph. However, these are not the only choices for coset representatives; in each case, there are two choices. It turns out that if you choose the other coset representative for each vertex of the cycle $C'$, you get a problem. These new coset representatives are $\theta^i\alpha^2\theta^{-1}$ for $i=0,1,2,3,4$. The two distinct vertices $u=\theta^{-1}H$ and $v=\alpha^2\theta^{-1}H$ give $u\theta=v\theta$ so these edges determined by $\theta$ can't both be a part of a 1-factor.

**Tree-like and Neighborhood Preserving Spanning Factorizations**.

We need a few more definitions to discuss spanning factorizations.

*Definitions*. Given words $\omega$ and $\mu$, we say $\mu\prec\omega$ if there exists a word $\alpha$ such that $\mu\alpha=\omega$. A spanning factorization $W=\{\omega_0=\varnothing,\omega_1=F_1,\omega_2=F_2,\ldots,\omega_d=F_d,\omega_{d+1},\ldots,\omega_{n-1}\}$ is *tree-like* if for every $\omega_j\in W$ and $\mu\prec\omega_j$ then there exists $i\le j$ such that $\mu=\omega_i$.

Given $\alpha\in W$ then by the definition of spanning factorization, for each $\beta\in W$ there exists a unique $\gamma\in W$ such that $u_0\gamma=(u_0\alpha)\beta$. We denote this $\gamma$ by $\alpha*\beta$. A spanning factorization is *neighborhood preserving* if for every $\alpha,\beta,F_i\in W$, we have $(\alpha*\beta)*F_i=\alpha*(\beta*F_k)$ for some $k$ (note that $k$ can depend on $\alpha,\beta$ and $F_i$).

**Properties of Digraphs with a Spanning 1-Factorization**

*Notation*. We assign acronyms to each of the relevant properties. We say a connected digraph is *VT* if it is vertex transitive; it is *SF* (*S1F)* if it has a spanning (1-)factorization; it is *TL* if it is SF and tree-like and it is *NP* if it is SF and neighborhood preserving.

*Theorem 1. (VT $\Rightarrow$SF, TL, NP) A connected vertex transitive diagraph $G=(\Gamma,S,H)$ has a tree-like, spanning factorization which is neighborhood preserving.*

*Proof.* We let $u_0=H$ and $u_i=s_iH$ for $1\le i\le d$ be its neighbors in $G=(\Gamma,S,H)$. We are going to form a factorization $F=\{F_i\}$ of $G$ and we start by assigning the edge $(H,s_iH)$ to the set $F_i$. We will label the vertices of $G$ by a breadth first lexicographically ordered tree of words $W(S)=\{\omega_i(S)\}$ in $S$ all of which are automorphisms. The elements of this tree have the property that the words are a set of unique coset representatives of $H$ in $G$. Note that paths in the tree are parsed from left to right because edges in the graph have the form $(\omega H,\omega s_iH)$ for $\omega\in W$. We know this is possible because we showed that in the connected case, above, $\Gamma=\langle S\rangle H$. $W(S)$ has these properties

1) every vertex $u_k$ in $G$ has the unique expression $u_k=\omega_k(S)(u_0)$;
2) $\omega_i\prec\omega_j\in W(S)$ implies $\omega_i\in W(S)$.

We have defined $u_0F_i = u_i = s_iH$ for $1 \le i \le d$. We expand this to include the edge $(\omega_k(S)(u_0), \omega_k(S)s_i(u_0)) = (u_0\omega_k(F), u_0\omega_k(F)F_i)$ in $F_i$ for every $\omega_k \in W(S)$. This gives us a set of factors $\{F_i\}$.

*Step 1*. First, we show that the factors are disjoint. Suppose that

$$(u_0\omega_l(F), u_0\omega_l(F)F_i) = (u_0\omega_k(F), u_0\omega_k(F)F_j).$$

Then $\omega_l(S)H = \omega_k(S)H$ so $l = k$. Then $\omega_k(S)s_iH = \omega_k(S)s_jH$ and because $\omega_k$ are automorphisms, $u_i = s_iH = s_jH = u_j$. This implies $i = j$ and so each edge is assigned to exactly one $F_i$.

*Step 2*. Next, we show that $u_iW(F) = u_0\omega_i(F)W(F)$ is a tree. If $u_0\omega_i(F)\alpha(F) = u_0\omega_i(F)\beta(F)$ for $\alpha(F), \beta(F), \omega_i(F) \in W(F)$, then $\omega_i(S)\alpha(S)H = \omega_i(S)\beta(S)H$. Since $\omega_i$ is an automorphism, $\alpha(S)H = \beta(S)H$. But since $\alpha$ and $\beta$ are unique coset representatives, this means $\alpha = \beta$.

*Step 3*. Finally, we show that this factorization is neighborhood preserving. Let $\alpha, \beta, F_i \in W$ and consider $(\alpha * \beta) * F_i$. There is $h \in H$ so that $\alpha * \beta = \alpha\beta h \in W$. Then there is $k \in H$ so that $\alpha\beta h * F_i = \alpha\beta h s_i k$. Note that $hs_ik \in HSH = SH$ so there are $l$ and $s_j$ such that $hs_ik \in s_jH$ which yields $\alpha\beta hs_ik \in \alpha\beta s_jH$. There exists $a$ such that $\beta s_jaH = \beta s_jH$ and $\beta s_ja \in W$ and there exists $b$ so that $\alpha\beta s_jaH = \alpha\beta s_jabH$ and $\alpha\beta s_jabH \in W$. Thus $\alpha\beta s_jabH = \alpha\beta hs_ikH$. This proves $(\alpha * \beta) * F_i = \alpha * (\beta * F_j)$.

*Theorem 2. (SF, NP $\Rightarrow$ VT) A digraph $G$ with a neighborhood preserving spanning factorization $W = \{\omega_0 = \varnothing, \omega_1 = F_1, \omega_2 = F_2, \ldots, \omega_d = F_d, \omega_{d+1}, \ldots, \omega_{n-1}\}$ is vertex transitive.*

*Proof.* We reverse the process in Theorem 2 but we do not need to use TL. For each $\omega \in W$, we define a function $\theta_\omega$ so that $\theta_\omega(u_0\omega_i) = u_0\omega * \omega_i$. To show that $\theta_\omega$ is an automorphism we have to show that it is a permutation that maps edges to edges. Suppose

$$u_0\omega * \omega_i = \theta_\omega(u_0\omega_i) = \theta_\omega(u_0\omega_j) = u_0\omega * \omega_j.$$

But then by the definition of spanning factorization $\omega_i = \omega_j$ so $\theta_\omega$ is a permutation. A typical edge has the form $(u_0\omega_i, u_0\omega_i * F_j)$ because $u_0\omega_i$ is a typical vertex and then the $(u_0\omega_i)F_j$ are its neighbors. Then $\theta_\omega(u_0\omega_i, u_0\omega_i * F_j) = (u_0\omega * \omega_i, u_0\omega * (\omega_i * F_j))$ and because $W$ is neighborhood preserving, $u_0\omega * (\omega_i * F_j) = u_0(\omega * \omega_i) * F_k$ for some $1 \le k \le d$. Thus $\theta_\omega$ maps the edge $(u_0\omega_i, u_0\omega_i * F_j)$ to the edge $(u_0\omega * \omega_i, u_0(\omega * \omega_i) * F_j)$. This shows that given any vertex $u_0\omega$ there is an automorphism $\theta_\omega$ which maps $u_0$ to $u_0\omega$ so $G$ is vertex transitive.

Putting these two theorems together gives us a characterization of vertex transitive graphs in terms of graphs with a spanning factorization. (This theorem is anticipated by a result of Mwambene in a 2006 paper [6] using the language of groupoids although this equivalence is not obvious.)

*Theorem 3. (SF, NP $\Leftrightarrow$ VT) A digraph $G$ is vertex transitive if and only if it has a neighborhood preserving spanning factorization* .

*Example 1. (S1F $\not\Rightarrow$ VT)* There exists a digraph that has a spanning 1-factorization but is not vertex transitive [7]. The graph has vertices 0 through 7 modulo 8. The degree is 2. The edges are the cycle $(i,i+1)$, the four edges $(i,i+2)$ for odd $i$ and the four edges $(i,i+4)$ for even $i$. The 1-factor $F_1$ is the cycle. The remaining edges are the 1-factor $F_2$. The spanning 1-factorization is then

$$W = \{\varnothing, F_1, F_2, F_1^3, F_1F_2F_1^7, F_1^5, F_1^6, F_1^7\}.$$

*Example 2. (VT $\not\Rightarrow$ S1F+TL)* There exists vertex transitive digraph of degree 2 which has no tree-like spanning 1-factorization. This digraph is $G=(\Gamma,S,H)$ with $\Gamma = A_5$, $S=\{(0,1,2,3,4),(0,3,4,2,1)\}$ and $H=\langle (0,2)(1,3)\rangle$. For details see [8].

*Open questions*. The digraph in Example 1 is not tree-like. We can't rule out

*VT $\Rightarrow$ S1F*
*SF, TL $\Rightarrow$ VT*
*S1F, TL $\Rightarrow$ VT*